\documentclass[10pt, reqno]{amsart}
\usepackage{graphicx, amssymb, amsmath, amsthm}
\numberwithin{equation}{section}
\usepackage{mathrsfs}
\usepackage{mathtools}
\usepackage{microtype}
\usepackage[T1]{fontenc}
\usepackage{lmodern}
\usepackage[colorlinks=true,linkcolor=blue,citecolor=blue,urlcolor=blue]{hyperref}
\usepackage[nameinlink,noabbrev]{cleveref}

\newtheorem{theorem}{Theorem}[section]
\newtheorem{proposition}[theorem]{Proposition}
\newtheorem{lemma}[theorem]{Lemma}

\theoremstyle{definition}

\theoremstyle{remark}
\newtheorem{remark}[theorem]{Remark}

\newcommand{\R}{\mathbb R}

\newcommand{\Sph}{\mathbb S^2}

\newcommand{\cF}{\mathcal F}
\newcommand{\cT}{\mathcal T}
\newcommand{\cW}{\mathcal W}

\newcommand{\dd}{\,d}

\newcommand{\norm}[1]{\left\|#1\right\|}

\newcommand{\ip}[2]{\left\langle#1,#2\right\rangle}

\DeclareMathOperator{\sech}{sech}
\DeclareMathOperator{\Dom}{D}
\DeclareMathOperator{\Rea}{Re}
\DeclareMathOperator{\Ima}{Im}

\title[Intertwining operators for relativistic Hardy operators]{Intertwining operators for relativistic Hardy operators and applications}

\author{YANSEN XIONG}
\address{Yansen Xiong:
Department of Mathematics, Beijing Institute of Technology, Beijing 100081}
\email{xiong\_yansen@bit.edu.cn}

\author{XIAOHUA YAO}
\address{Xiaohua Yao:
Department of Mathematics and Hubei Province Key Laboratory of Mathematical Physics, Central China Normal University, Wuhan  430079}
\email{yaoxiaohua@mail.ccnu.edu.cn}

\author{JUNYONG ZHANG}
\address{Junyong Zhang:
Department of Mathematics, Beijing Institute of Technology, Beijing 100081}
\email{zhang\_junyong@bit.edu.cn}


\begin{document}

\begin{abstract}
We study the scaling-critical relativistic Hardy operator
$H_a=\sqrt{-\Delta}+a|x|^{-1}$ on $\mathbb{R}^3$ for $a>-2/\pi$. We construct an explicit unitary intertwiner between
$H_a$ and $\sqrt{-\Delta}$, obtained from a holomorphic Mellin
factorization in each angular channel.
We prove estimates for the
factor that are uniform simultaneously in the Mellin variable and the
angular momentum; they yield $L^p$ bounds for the intertwiner and its
inverse on an explicit coupling-dependent interval, as well as
Sobolev isomorphisms.   As an application, we establish Strichartz
estimates for the massless Coulomb evolution, thereby resolving the
three-dimensional scalar first-order case left open by
Mizutani--Yao \cite{MY}. For attractive couplings, the spatial
restriction $p<3/\delta_a$ is necessary, where $\delta_a\in(0,1)$ is
the singular ground-state exponent. \end{abstract}

\maketitle

\section{Introduction}

\subsection{Background and motivation}

The relativistic Schr\"odinger operator with a Coulomb potential
provides a scalar model for a relativistic particle in the field of a
point charge. In units in which the relevant physical constants are one,
its kinetic energy is $\sqrt{-\Delta+m^2}$, with mass $m\geq0$, while
the Coulomb interaction has the form $a|x|^{-1}$. For $a<0$ the
interaction is attractive, and for $a>0$ it is repulsive. The spectral
theory of the massive operator was studied by Herbst \cite{Herbst}.
In this paper we consider its massless counterpart
\[
 H_a=|D|+a|x|^{-1},\qquad |D|=\sqrt{-\Delta},
 \qquad x\in\R^3.
\]
The first-order kinetic energy and the Coulomb potential have the same
homogeneity: under spatial dilation both scale like inverse length.
Consequently, neither term is a lower-order perturbation of the other
at any scale. This scaling-critical balance, present at every coupling
$a$, is central to both the physical model and the analysis of the
propagator $e^{-itH_a}$.

There is also a natural comparison with the Dirac--Coulomb operator
\[
 \mathcal D_{m,a}=-i\boldsymbol\alpha\cdot\nabla+m\beta
                    +a|x|^{-1}I_4,
\]
acting on four-component spinors, where $\boldsymbol\alpha$ and
$\beta$ are Dirac matrices. The free Dirac operator satisfies
$\mathcal D_{m,0}^{\,2}=(-\Delta+m^2)I_4$, and hence
$|\mathcal D_{m,0}|=\sqrt{-\Delta+m^2}\,I_4$. In particular,
$|\mathcal D_{0,0}|=|D|I_4$, so the scalar kinetic energy is the
absolute value of the free Dirac kinetic energy. The two interacting
models are nevertheless different: the Dirac operator retains its
spinor structure and positive and negative spectral branches, and the
Coulomb potential does not in general commute with the free spectral
projections. Thus $H_a$ is not being identified with a restriction or
absolute value of $\mathcal D_{0,a}$, and our scalar estimates do not
by themselves imply estimates for Dirac--Coulomb evolution. The latter
problem is also highlighted in \cite[Remark~1.14]{MY}.

The scaling-critical nature of $H_a$ should be distinguished from
criticality of the coupling constant. The sharp Hardy inequality in
three dimensions gives
\begin{equation}\label{hardy_ineq}
    \frac{2}{\pi}\int_{\R^3}\frac{|f(x)|^2}{|x|}\,\dd x
 \leq \bigl\||D|^{1/2}f\bigr\|_2^2,
 \qquad f\in C_c^\infty(\R^3),
\end{equation}
and identifies $a=-2/\pi$ as the critical attractive coupling for the
scalar quadratic form; see \cite{Herbst}. We work throughout with
$a>-2/\pi$, where the form is comparable to that of $|D|$. Thus our
results treat the scaling-critical massless model in the
\emph{subcritical coupling} regime, not the endpoint $a=-2/\pi$.
It follows from Hardy's inequality \eqref{hardy_ineq} that the relativistic Hardy operators $H_a$ can be realized as self-adjoint operators on $L^{2}(\R^3)$ (see Section \ref{sec:mellin} for the precise definition).
More generally, the fractional Hardy family
\[
 (-\Delta)^\sigma+a|x|^{-2\sigma},\qquad 0<\sigma<d/2,
\]
has the same balance between kinetic and potential energies. The
relativistic Coulomb operator is its three-dimensional member with
$\sigma=1/2$.

The analysis of dispersive equations with singular critical potentials
includes the work of Burq, Planchon, Stalker, and Tahvildar-Zadeh
\cite{BPST} on inverse-square potentials. In the fractional Hardy
family above, the local second-order member $\sigma=1$ is
\[
 \mathcal L_a=-\Delta+a|x|^{-2}.
\]
Miao, Su, and Zheng \cite{MSZ} constructed stationary wave operators
for $\mathcal L_a$ and established their $W^{s,p}$-boundedness, with
applications to dispersive and smoothing estimates. Their construction
uses the explicit Bessel-function spectral representation of the local
radial Schr\"odinger operator. It therefore cannot be applied as it
stands to the nonlocal first-order operator $H_a$, for which there is no
corresponding local radial differential equation.

For the fractional Hardy family, Mizutani and Yao \cite{MY} established resolvent and Kato
smoothing estimates for $0<\sigma<d/2$ and proved Strichartz estimates
for $\sigma>1/2$ under subcritical coupling. They explicitly ask whether
Strichartz estimates hold in the first-order case for the
relativistic Coulomb operator, even for small nonzero coupling, and more
generally for $0<\sigma\leq1/2$; see \cite[Remarks~1.14 and 1.16]{MY}.
Theorem~\ref{thm:main} addresses the three-dimensional scalar
$\sigma=1/2$ instance of that question for all $a>-2/\pi$, in the
admissible range and with the spatial restriction stated below. It does
not settle the cases $0<\sigma<1/2$, the critical coupling, or the Dirac
problem.

The first-order case is not accessible by a direct application of the
perturbative argument in \cite{MY}. Indeed, its smoothing estimate
specializes in three dimensions to
\[
 \bigl\||x|^{-1/2+\gamma}|D|^\gamma
 e^{-itH_a}f\bigr\|_{L_t^2L_x^2}
 \lesssim \|f\|_{L_x^2},\qquad -1<\gamma<0,
\]
but does not include $\gamma=0$, the value used in the perturbative
Strichartz argument for $\sigma>1/2$.

More recently, Fanelli, Su, Wang, Zhang, and Zheng \cite{FSWZZ}
developed spectrally projected intertwining operators for a broad class
of second-order scaling-critical Schr\"odinger Hamiltonians, extending
the method beyond constant inverse-square perturbations to nonconstant
spherical electromagnetic potentials. This framework further
demonstrates the usefulness of intertwining methods, but its local
second-order spectral representation does not provide an intertwiner
for the nonlocal kinetic energy $|D|$.

Our approach is nonperturbative. Its main new ingredient is an explicit
holomorphic factorization of the channelwise Mellin-symbol ratio
$(M_L+a)/M_L$. The factor solves a vertical difference relation and is
controlled uniformly in both the Mellin variable and the angular
parameter $L$. This simultaneous control is essential: after the
$L^p$ polar change of variables, the holomorphic strip becomes exactly
the interval of admissible $p$, while the same strip yields the
Sobolev isomorphisms. The resulting unitary conjugation transfers the
free half-wave estimates to the perturbed evolution without invoking
the missing endpoint smoothing estimate.

The singular behavior in the attractive case also constrains the range
of the result. Frank, Lieb, and Seiringer \cite{FLS} developed a
ground-state representation for fractional Hardy forms, and Bogdan,
Grzywny, Jakubowski, and Pilarczyk \cite{BGJP} obtained sharp heat-kernel
bounds. For $-2/\pi<a<0$, the corresponding ground-state profile is
$|x|^{-\delta_a}$, where $\delta_a\in(0,1)$ is defined below. Its local
integrability requires $p\delta_a<3$; the heat-kernel lower bound shows
that this restriction is necessary for the Strichartz estimates proved
here. Related angular-channel ground-state representations and heat
kernel estimates were obtained by Bogdan and Merz
\cite{BogdanMerz2024,BogdanMerz2025}. For the Sobolev comparison, we
also use results of Frank, Merz, and Siedentop \cite{FMS} and Merz
\cite{Merz}. Our multiplier construction supplements these results with
the uniform angular estimates needed for the oscillatory evolution.

\subsection{Main results}
Throughout this paper, we assume $a>-2/\pi$. To describe the restriction on the
spatial exponent in the attractive case, define
\begin{equation}\label{eq:delta}
 \delta_a=
 \begin{cases}
 0, & a\geq0,\\
  \delta, & 0>a>-2/\pi,
 \end{cases}
\end{equation}
where $\delta\in(0,1)$ satisfies $a=-c(\delta)$, with
\begin{equation}\label{eq:c}
 c(\delta)
 =2\frac{
 \Gamma\bigl((\delta+1)/2\bigr)
 \Gamma\bigl((3-\delta)/2\bigr)
 }{
 \Gamma(\delta/2)
 \Gamma\bigl((2-\delta)/2\bigr)
 }
\end{equation}
The function $c$ is expressed in terms of the Gamma function $\Gamma$.
The function $c$ increases strictly on $(0,1)$ from $0$ to 
$2/\pi$.
Thus, for every
$-2/\pi<a<0$, the parameter $\delta_a$ is uniquely determined.
For $a\geq0$, the convention $\delta_a=0$ is adopted solely to
unify the statements of the estimates below.

We also introduce the interval
\begin{equation}\label{eq:interval}
 \begin{aligned}
 I_a
 &:=\left\{p\in(1,\infty):
 \left|\frac32-\frac3p\right|<\frac32-\delta_a\right\}
 =
 \begin{cases}
 (1,\infty),&a\geq0,\\[2pt]
 \bigl(3/(3-\delta_a),\,3/\delta_a\bigr),&a<0.
 \end{cases}
 \end{aligned}
\end{equation}

We realize $H_a$ as the nonnegative self-adjoint
operator associated with the closed quadratic form
\[
 q_a[f]=\||D|^{1/2}f\|_2^2
 +a\int_{\R^3}|x|^{-1}|f(x)|^2\,\dd x,
 \qquad \Dom(q_a)=H^{1/2}(\R^3).
\]
To state our main result, we define the intertwiner using logarithmic
polar coordinates. Let $\cT:L^2(\R^3)\to L^2(\R\times\Sph)$ be the
unitary transform
\begin{equation}\label{eq:intro-polar}
 (\cT f)(y,\omega)=e^{3y/2}f(e^y\omega).
\end{equation}
Fix an orthonormal basis $\{Y_{\ell m}:\ell\geq0,\ -\ell\leq m\leq\ell\}$
of spherical harmonics on $\Sph$, with $Y_{\ell m}$ of degree $\ell$.
For $L\geq1$, set
\begin{equation}\label{ML}
 M_L(z)=2\frac{\Gamma((L+1+iz)/2)\Gamma((L+1-iz)/2)}
 {\Gamma((L+iz)/2)\Gamma((L-iz)/2)},\end{equation}
 and
 \begin{equation}\label{bL}
 b_L(\xi)=\log\!\left(1+\frac a{M_L(\xi)}\right),
\end{equation}
where the logarithm is real on $\R$. Define
\begin{equation}\label{WL}
 K(x)=\frac\pi2\sech^2(\pi x),\quad
 w_L(\xi)=\exp\!\left(-i\int_0^\xi (K*b_L)(s)\,ds\right).
\end{equation}
The functions $w_L$ are unimodular on $\R$; their analytic construction
and further properties are established in Lemma~\ref{lem:factor},
see \eqref{eq:construction}.

For $u\in L^2(\R\times\Sph)$, write its spherical harmonic expansion as
\begin{equation}\label{eq:intro-expansion}
 u(y,\omega)=\sum_{\ell=0}^\infty\sum_{m=-\ell}^{\ell}
 u_{\ell m}(y)Y_{\ell m}(\omega),\qquad
 u_{\ell m}(y)=\int_{\Sph}u(y,\omega)\overline{Y_{\ell m}(\omega)}\,d\omega,
\end{equation}
where the expansion converges in $L^2(\R\times\Sph)$. With
$D_y=-i\partial_y$, define $w_{\ell+1}(D_y)$ as the Fourier multiplier
in the $y$ variable with symbol $w_{\ell+1}$:
\begin{equation}\label{eq:intro-multiplier}
 \cF_y\!\left[w_{\ell+1}(D_y)u_{\ell m}\right](\xi)
 =w_{\ell+1}(\xi)\,\cF_yu_{\ell m}(\xi).
\end{equation}
We then define
\begin{equation}\label{eq:intro-cylinder-W}
 (\cW_a u)(y,\omega)=\sum_{\ell=0}^\infty\sum_{m=-\ell}^{\ell}
 \bigl(w_{\ell+1}(D_y)u_{\ell m}\bigr)(y)Y_{\ell m}(\omega),
\end{equation}
with convergence in $L^2(\R\times\Sph)$. By Plancherel's theorem,
the orthogonality of the spherical harmonics and $|w_L|=1$ on $\R$,
$\cW_a$ is unitary; its inverse is given by the same formula with
$w_{\ell+1}$ replaced by $\overline{w_{\ell+1}}$. Finally, set
\begin{equation}\label{eq:W}
 W_a=\cT^{-1}\cW_a\cT.
\end{equation}
Thus $W_a$ is unitary on $L^2(\R^3)$. Proposition~\ref{prop:intertwiner}
establishes its intertwining identity \eqref{eq:structure}. This is a
particular intertwiner; we do not identify it with a time-dependent
scattering wave operator.

Our first main result is the following. 
\begin{theorem}[$L^p$ bounds for intertwining operators]\label{thm:intertwining}
Let $a>-2/\pi$, and let $I_a$ be given by \eqref{eq:interval}. The
operator $W_a$ in \eqref{eq:W} satisfies
$H_a=W_a |D| W_a^*$.
For every $p\in I_a$, both $W_a$
and $W_a^*$ extend to bounded operators on $L^p(\R^3)$, and these
extensions are mutually inverse. In particular, for every Borel function
$m$, $m(H_a)=W_a m(|D|)W_a^*$ on the corresponding spectral domain.
\end{theorem}
\begin{remark}
    We show that
the function $w_L$ satisfies uniform estimates in both
the variable $\xi$ and the angular parameter $L$. These estimates
allow us to apply the joint spectral multiplier theorem of
Wr\'obel~\cite{Wrobel}, which implies $L^p$-boundedness of the intertwining operators $W_a$.
\end{remark}

As an application, we transfer the free half-wave estimates to $H_a$.
Unless otherwise indicated, spatial norms are over $\R^3$, time norms
are over $\R$, and $L_t^rL_x^p$ denotes $L^r(\R_t;L^p(\R_x^3))$.
We interpret $3/\delta_a=\infty$ when $\delta_a=0$.

\begin{theorem}[Strichartz estimates]
\label{thm:main}
Let $a>-2/\pi$, and suppose that
\begin{equation}\label{eq:admissible}
 2<r\leq\infty,\qquad 2\leq p< 3/\delta_a,\qquad
 1/r+1/p\leq1/2,
\end{equation}
and the scaling condition
\begin{equation}
 s=\frac32-\frac1r-\frac3p.
\end{equation}
Then
\begin{equation}\label{eq:main}
 \norm{e^{-itH_a}f}_{L_t^rL_x^p}
 \leq C_{a,r,p}\norm{f}_{\dot H^s(\R^3)}
\end{equation}
for every $f\in\dot H^s(\R^3)$. 
\end{theorem}

\begin{remark}
The restriction $p< 3/\delta_a$ is necessary. More precisely, if
$-2/\pi<a<0$ and $p\geq 3/\delta_a$, then no estimate of the form
\eqref{eq:main} holds.
\end{remark}

The proof begins with an explicit factorization adapted to the
Coulomb singularity. In logarithmic polar coordinates, the quadratic
form of $|D|$ in the angular channel $\ell$ has Mellin symbol
$M_{\ell+1}$, whereas the potential adds $a$; see
Lemma~\ref{lem:mellinform}. The factor $w_L$ in
\eqref{eq:intro-multiplier} is constructed from the logarithm of the
ratio $(M_L+a)/M_L$ by convolution with $K$ and an exponential
primitive. Its boundary identity \eqref{eq:factoridentity} yields the
exact conjugation \eqref{eq:structure} at the level of quadratic forms.

A further key ingredient is control of $w_L$ and its inverse uniformly
in both the Mellin variable and the angular parameter, on the strip
$|\Ima z|<3/2-\delta_a$; see Lemma~\ref{lem:factor}. The $L^p$ polar
isometry shifts the Mellin variable to
$\xi-i(3/2-3/p)$, so the available strip gives precisely the range
$p\in I_a$ in \eqref{eq:interval}. The mixed estimates permit the use
of Wr\'obel's joint multiplier theorem \cite{Wrobel} to sum the angular
channels, proving Theorem~\ref{thm:intertwining}. The same strip
estimates also yield the Sobolev comparison
\eqref{eq:two-sided-sobolev}. The intertwining identity then transfers
the free half-wave estimates to prove Theorem~\ref{thm:main}; the
necessity of its attractive spatial restriction follows from the heat
kernel lower bound in Proposition~\ref{prop:attractive-obstruction}.

\section{Quadratic forms and intertwining operators}\label{sec:mellin}

In this section we construct the unitary intertwiner introduced in
\eqref{eq:W} and prove the conjugation identity
\eqref{eq:structure} below. The construction starts with a Mellin
representation of the quadratic form, followed by an analytic
factorization in each angular channel.

\subsection{Quadratic forms and Mellin representation}

We recall the quadratic form defining $H_a$ in the introduction:
\begin{equation}\label{eq:form}
 q_a[f]
 =\norm{|D|^{1/2}f}_2^2
  +a\int_{\R^3}|x|^{-1}|f(x)|^2\,\dd x,
 \qquad
 \Dom(q_a)=H^{1/2}(\R^3).
\end{equation}
The sharp Hardy inequality in three dimensions is
\begin{equation}\label{eq:hardy}
 \frac{2}{\pi}\int_{\R^3}|x|^{-1}|f(x)|^2\,\dd x
 \leq \int_{\R^3}\bigl||D|^{1/2}f(x)\bigr|^2\,\dd x.
\end{equation}
For $a>-2/\pi$, combining \eqref{eq:hardy} with \eqref{eq:form} yields
\begin{equation}\label{eq:formcomp}
 c_a q_0[f]\leq q_a[f]\leq C_a q_0[f],
\end{equation}
where $c_a=\min\{1,1+\pi a/2\}>0$ and
$C_a=\max\{1,1+\pi a/2\}$.
In particular, $q_a$ is closed on
$H^{1/2}(\R^3)$ and determines the Friedrichs realization of
$H_a$ by \cite{Kato}.

The dilation
\[
 (U_\lambda f)(x)=\lambda^{3/2}f(\lambda x),
 \qquad \lambda>0,
\]
satisfies
\begin{equation}\label{eq:scaling}
 U_\lambda^{-1}H_aU_\lambda=\lambda H_a,
 \qquad
 e^{-itH_a}U_\lambda=U_\lambda e^{-i\lambda tH_a}.
\end{equation}
Indeed, both terms in $q_a$ scale by the same factor $\lambda$.
The corresponding operator identity follows from the representation
theorem for closed quadratic forms. 

We use the orthonormal basis of  spherical harmonics $\{Y_{\ell m}\}$ and the unitary
logarithmic polar transform $\cT$ from \eqref{eq:intro-polar}:
\begin{equation}\label{eq:polar}
 (\cT f)(y,\omega)=e^{3y/2}f(e^y\omega),
 \qquad
 \cT:L^2(\R^3)\longrightarrow L^2(\R\times\Sph).
\end{equation}
Let $D_y=-i\partial_y$, and let $\cF_y$ denote the unitary Fourier
transform in $y$, unless another normalization is specified. Define
the angular operator
\[
 \Lambda=\sqrt{-\Delta_{\Sph}+\frac14}+\frac12,
\]
so that
\[
 \Lambda Y_{\ell m}=(\ell+1)Y_{\ell m}.
\]
Recall the Mellin symbol from \eqref{ML}:
\begin{equation}\label{eq:M}
 M_L(z)=
 2\frac{
 \Gamma\bigl((L+1+iz)/2\bigr)
 \Gamma\bigl((L+1-iz)/2\bigr)
 }{
 \Gamma\bigl((L+iz)/2\bigr)
 \Gamma\bigl((L-iz)/2\bigr)
 }.
\end{equation}

\begin{lemma}\label{lem:mellinform}
Let $f\in\Dom(q_a)$ and write, as in \eqref{eq:intro-expansion},
\[
u=\mathcal{T}f=\sum_{\ell,m}u_{\ell m}(y)Y_{\ell m}(\omega).
\]
Then
\begin{equation}\label{eq:mellinform}
 q_a[f]=\sum_{\ell,m}\int_\R
 (M_{\ell+1}(\xi)+a)
 \left|\cF_y\big(e^{-y/2}u_{\ell m}(y)\big)(\xi)\right|^2\dd\xi.
\end{equation}
Moreover, smooth functions supported in compact spatial annuli and
containing only finitely many angular modes form a common core for
$q_0$ and $q_a$.
\end{lemma}

\begin{proof}

We first verify the representation for smooth functions supported in a
compact spatial annulus and involving only finitely many angular modes.
We then show that this class is a common form core.

We use the classical Fourier transform formula for homogeneous
spherical harmonics \cite{SteinWeissBook}:
\begin{equation}\label{eq:homogeneous}
 \cF\bigl(|x|^{-\beta}Y_{\ell m}(x/|x|)\bigr)(\xi)
 =(-i)^\ell 2^{3/2-\beta}
 \frac{\Gamma((\ell+3-\beta)/2)}
      {\Gamma((\ell+\beta)/2)}
 |\xi|^{\beta-3}Y_{\ell m}(\xi/|\xi|),
\end{equation}
which holds, as an identity of tempered distributions, for
\(0<\Rea\beta<3\). It follows that 
\[
 \cF\bigl(|D|\bigl(|x|^{-\beta}Y_{\ell m}(x/|x|)\bigr)\bigr)(\xi)
 =(-i)^\ell 2^{3/2-\beta}
 \frac{\Gamma((\ell+3-\beta)/2)}
      {\Gamma((\ell+\beta)/2)}
 |\xi|^{\beta-2}Y_{\ell m}(\xi/|\xi|)
\]
and
\[
 \cF\bigl(|x|^{-\beta-1}Y_{\ell m}(x/|x|)\bigr)(\xi)
 =(-i)^\ell 2^{1/2-\beta}
 \frac{\Gamma((\ell+2-\beta)/2)}
      {\Gamma((\ell+\beta+1)/2)}
 |\xi|^{\beta-2}Y_{\ell m}(\xi/|\xi|),
\]
Comparison with \eqref{eq:homogeneous} gives
\begin{equation}\label{eq:homaction}
 |D|\bigl(|x|^{-\beta}Y_{\ell m}\bigr)
 =
 2\frac{\Gamma((\ell+\beta+1)/2)
         \Gamma((\ell+3-\beta)/2)}
        {\Gamma((\ell+\beta)/2)
         \Gamma((\ell+2-\beta)/2)}
 |x|^{-\beta-1}Y_{\ell m}
\end{equation}
for $0<\Rea\beta<2$.

Now we fix an angular mode \((\ell,m)\) and we write
\[
 f(r\omega)=r^{-1}h(\log r)Y_{\ell m}(\omega).
\]
Fourier inversion in \(y\) gives
\[
h(y)=\frac1{\sqrt{2\pi}}\int_{\mathbb R}
e^{iy\xi}\widehat h(\xi)\,d\xi,
\]
which yields
\begin{equation}\label{eq:fromega}
f(r\omega)=\frac1{\sqrt{2\pi}}\int_{\mathbb R}
\widehat h(\xi)\,r^{-1+i\xi}Y_{\ell m}(\omega)\,d\xi.
\end{equation}
It follows from \eqref{eq:homaction} and \eqref{eq:M} that 
\[
|D|\bigl(r^{-1+i\xi}Y_{\ell m}\bigr)
=M_{\ell+1}(\xi)\,r^{-2+i\xi}Y_{\ell m}.
\]
Applying this identity to \eqref{eq:fromega} and using \eqref{eq:M}
yields
\[
|D|f(r\omega)=r^{-2}k(\log r)Y_{\ell m}(\omega),
\qquad
\widehat k(\xi)=M_{\ell+1}(\xi)\widehat h(\xi),
\]
which implies
\[
\begin{aligned}
q_0[f]
&=\int_{\mathbb R^3}|D|f(x)\overline{f(x)}\,dx\\
&=\int_0^\infty\!\int_{\mathbb S^2}
\bigl[r^{-2}k(\log r)Y_{\ell m}\bigr]
\bigl[r^{-1}\overline{h(\log r)}\,\overline{Y_{\ell m}}\bigr]
r^2\,d\omega\,dr\\
&=\int_{\mathbb R}k(y)\overline{h(y)}\,dy=\int_{\mathbb R}\widehat k(\xi)\overline{\widehat h(\xi)}\,d\xi\\
&=\int_{\mathbb R}M_{\ell+1}(\xi)|\widehat h(\xi)|^2\,d\xi.
\end{aligned}
\]
On the other hand, we have
\[
\begin{aligned}
a\int_{\mathbb R^3}|x|^{-1}|f(x)|^2\,dx
&=a\int_0^\infty r^{-1}r^{-2}|h(\log r)|^2r^2\,dr\\
&=a\int_{\mathbb R}|\widehat h(\xi)|^2\,d\xi.
\end{aligned}
\]
Consequently,
\[
q_a[f]=\int_{\mathbb R}
\bigl(M_{\ell+1}(\xi)+a\bigr)|\widehat h(\xi)|^2\,d\xi.
\]
Write
\[
u=\mathcal T f,\qquad u(y,\omega)=e^{3y/2}f(e^y\omega).
\]
In this channel, \(u_{\ell m}(y)=e^{y/2}h(y)\), so
\(h=e^{-y/2}u_{\ell m}\). Summing over the orthogonal spherical
harmonic channels gives
\[
q_a[f]=\sum_{\ell,m}\int_{\mathbb R}
\bigl(M_{\ell+1}(\xi)+a\bigr)
\left|\mathcal F_y\!\left(e^{-y/2}u_{\ell m}\right)(\xi)\right|^2\,d\xi.
\]
This proves \eqref{eq:mellinform} as desirable.

It remains to identify the maximal domain of the quadratic form on the
right-hand side. By the standard asymptotics for ratios of Gamma
functions \cite{DLMF},
\begin{equation}\label{eq:Mreal}
 M_L(\xi)\asymp L+|\xi|,
 \qquad L\geq1,\quad \xi\in\mathbb R.
\end{equation}
Hence the free cylinder form associated with \(M_\Lambda(D_y)\) has
domain
\[
 \left\{
 u\in L^2(\mathbb R\times\Sph):
 h=e^{-y/2}u,\quad
 \|\Lambda^{1/2}h\|_{L^2}^2
 +\||D_y|^{1/2}h\|_{L^2}^2<\infty
 \right\}.
\]
The equivalence \eqref{eq:Mreal} also shows that the corresponding
graph norm is equivalent to the natural Hilbert-valued
\(H^{1/2}(\mathbb R)\) norm with the angular weight \(\Lambda^{1/2}\).
In particular, multiplication by \(e^{-y/2}\), with its natural
maximal domain, is closed. Standard cutoff, angular truncation, and
mollification arguments therefore show that smooth functions compactly
supported in \(y\) and involving only finitely many angular modes form
a core for this maximal form.

Under the inverse transformation \(\mathcal T^{-1}\), compact support in
\(y\) is equivalent to support in a compact spatial annulus. Hence the
same class is a form core for the Euclidean free form \(q_0\). More
explicitly, for a fixed smooth function, truncation near \(x=0\) and
\(x=\infty\) can be performed with an error tending to zero in
\(\dot H^{1/2}(\mathbb R^3)\), while on a fixed annulus the spherical
harmonic partial sums converge in \(H^1\), and hence in \(H^{1/2}\).
Thus the two realizations of the free form have the same closure.

Finally, the Hardy form comparison \eqref{eq:formcomp} implies that
\(q_a\) and \(q_0\) have the same form domain. Since
\(M_L(\xi)+a\) is comparable to \(M_L(\xi)\), uniformly for
\(L\geq1\) and \(\xi\in\mathbb R\), the quadratic form on the right-hand
side of \eqref{eq:mellinform} has the same maximal domain. Therefore
the identity extends from the common core to the entire form domain
by closedness. This completes the proof.
\end{proof}

\subsection{The logarithmic symbol}

By \eqref{eq:mellinform}, the Coulomb potential replaces $M_L$ by
$M_L+a$ in each angular channel. Our construction factors this
perturbation through a holomorphic function satisfying the shifted
boundary identities in Lemma~\ref{lem:factor} below. We first control the
logarithm of the ratio, extending the real-valued symbol in
\eqref{bL} to a complex strip. Uniformity in the continuous parameter
$L$ will be essential when the angular channels are recombined in
Section~\ref{sec:multipliers}.

\begin{lemma}\label{lem:logsymbol}
For each $L\geq1$, let $b_L$ be the real-valued function in
\eqref{bL}.
Then $b_L$ admits a single-valued holomorphic branch on the strip
\[
\{z\in\mathbb C:|\Ima z|< 1-\delta_a\},
\]
which is real-valued and even on the real axis. Furthermore, write
$\zeta=x+iy$. For every $\eta\in[0,1-\delta_a)$ and all nonnegative
integers $j,k$, there exists a constant $C_{a,\eta,j,k}>0$, independent
of $L$, $x$, and $y$, such that
\begin{equation}\label{eq:blogbounds}
 \left|\partial_\zeta^j\partial_L^k b_L(\zeta)\right|
 \leq C_{a,\eta,j,k}(L+|x|)^{-1-j-k},
\end{equation}
where $\zeta=x+iy$, $L\geq1$, $x\in\mathbb R$, and $|y|\leq\eta$.
\end{lemma}

\begin{proof}
We first record the symbol estimates for $M_L$. The precise
gamma-ratio expansion we use is \cite[Eq.~(5.11.13)]{DLMF}:
\begin{equation}\label{eq:gamma-ratio-asymptotic}
 \frac{\Gamma(w+\alpha)}{\Gamma(w+\beta)}
 \sim w^{\alpha-\beta}
 \sum_{n=0}^{\infty}\frac{G_n(\alpha,\beta)}{w^n},
 \qquad |w|\longrightarrow\infty,
\end{equation}
uniformly in every closed sector $|\arg w|\leq\pi-\varepsilon$.
For $\zeta=x+iy$, set
\[
 w_\pm=\frac{L\pm i\zeta}{2}
       =\frac{L\mp y\pm ix}{2}.
\]
If $L\geq1$ and $|y|\leq\eta<1$, then
\[
 \Rea w_\pm\geq\frac{1-\eta}{2},
 \qquad |w_\pm|\asymp_\eta L+|x|.
\]
Since
\[
 M_L(\zeta)
 =2\frac{\Gamma(w_++1/2)}{\Gamma(w_+)}
    \frac{\Gamma(w_-+1/2)}{\Gamma(w_-)},
\]
we apply \eqref{eq:gamma-ratio-asymptotic} with
$(\alpha,\beta)=(1/2,0)$ to both factors. Differentiating the resulting
uniform asymptotic expansions, or equivalently using Cauchy's formula
on slightly smaller sectors, gives, for all nonnegative integers $j,k$,
\begin{align}
 \bigl|\partial_\zeta^j\partial_L^k M_L(\zeta)\bigr|
 &\leq C_{\eta,j,k}(L+|x|)^{1-j-k}, \label{eq:MLsymbol}\\
 \bigl|\partial_\zeta^j\partial_L^k [M_L(\zeta)]^{-1}\bigr|
 &\leq C_{\eta,j,k}(L+|x|)^{-1-j-k}, \label{eq:MLinvsymbol}
\end{align}
for $\zeta=x+iy$, $L\geq1$, $x\in\mathbb R$, and $|y|\leq\eta$.
The reciprocal estimate follows from the same expansion, whose leading
term is nonzero. The remaining region in which $L+|x|$ is bounded is
compact. Since all Gamma-function arguments remain in a fixed right
half-plane, both estimates there follow from holomorphy and the fact
that the Gamma function has no zeros.

For $\Rea w>0$, define
\[
 R(w):=\frac{\Gamma(w+1/2)}{\Gamma(w)}.
\]
The beta integral and the binomial expansion of
$(1-t)^{-1/2}$ give
\begin{equation}\label{eq:Rseries}
 R(w)
 =\frac{w}{\sqrt\pi}\int_0^1
   t^{w-1/2}(1-t)^{-1/2}\,\dd t
 =\frac1{\sqrt\pi}\sum_{n=0}^{\infty}
 4^{-n}\binom{2n}{n}\frac{w}{w+n+1/2},
\end{equation}
which converges locally uniformly in the right half-plane. Writing $w=u+iv$ with $u>0$, we have
\[
 \Rea\frac{w}{w+n+1/2}
 =
 \frac{u(u+n+1/2)+v^2}
 {(u+n+1/2)^2+v^2}>0
\]
and
\[
 \Ima\frac{w}{w+n+1/2}
 =
 \frac{v(n+1/2)}
 {(u+n+1/2)^2+v^2}.
\]
Consequently, $\Rea R(w)>0$, and $\Ima R(w)$ has the same sign as
$\Ima w$.

We note that
\[
 M_L(z)
 =2R\left(\frac{L+iz}{2}\right)
    R\left(\frac{L-iz}{2}\right), \quad z=\zeta=x+iy.
\]
If $|\Ima z|<L$, then the two arguments of $R$ have positive real parts
and imaginary parts of opposite signs. Hence, their images under $R$ have positive real parts and imaginary parts of opposite signs, 
which implies $\Rea M_L(z)>0$.

Suppose first that $a>0$. Then $\delta_a=0$, and
$\Rea M_L(z)>0$ throughout $|\Ima z|<1-\delta_a$.
We see that
\[
 b_L(z)=\int_0^a\frac{\dd v}{M_L(z)+v}.
\]
The integrand is holomorphic in $z$, and
\[
 |M_L(z)+v|^2
 =|M_L(z)|^2+2v\Rea M_L(z)+v^2
 \geq |M_L(z)|^2
\]
for $v\in[0,a]$. It follows from \eqref{eq:MLinvsymbol} that
\[
 \bigl|\partial_z^j\partial_L^k
       (M_L(x+iy)+v)^{-1}\bigr|
 \leq C_{\eta,j,k}(L+|x|)^{-1-j-k},
\]
uniformly for $v\in[0,a]$. Differentiation under the integral sign
therefore yields \eqref{eq:blogbounds}. Moreover, the integral defines
a single-valued holomorphic branch of
$\log(1+a/M_L(z))$. Since $M_L$ is positive and even on the real axis,
$b_L(z)$ is real-valued and even there.

If $a=0$, then $b_L\equiv0$, and the conclusion is immediate.

Now assume that $a<0$. We note that
\begin{equation}\label{eq:betaM}
 \frac1{M_L(z)}
 =\frac2\pi\int_0^\infty\int_0^\infty
 \frac{e^{-L(v+w)}e^{-iz(v-w)}}
 {\sqrt{1-e^{-2v}}\sqrt{1-e^{-2w}}}
 \,\dd v\,\dd w,
 \qquad |\Ima z|<L,
\end{equation}
which implies
\begin{equation}\label{eq:reciprocalstrip}
 \bigl|M_L(x+iy)^{-1}\bigr|
 \leq M_L(iy)^{-1}
 \leq M_1(iy)^{-1}.
\end{equation}
Exchanging \(v\) and \(w\) in \eqref{eq:betaM} and averaging gives
\begin{equation}
 \frac1{M_L(iy)}
 =\frac2\pi\int_0^\infty\int_0^\infty
 \frac{e^{-L(v+w)}\cosh(y(v-w))}
 {\sqrt{1-e^{-2v}}\sqrt{1-e^{-2w}}}
 \,\dd v\,\dd w.
\end{equation}
Hence $M_L(iy)^{-1}$ is increasing in
$|y|$ for $|y|<L$. 
The gamma-function identities in \eqref{eq:M} and the definition
\eqref{eq:c} give
\begin{equation}\label{eq:Mone}
 M_1(z)=z\coth(\pi z/2),\qquad M_1(0)=2/\pi,
 \qquad M_1(i(1-\delta_a))=c(\delta_a)= -a.
\end{equation}
Applying \eqref{eq:reciprocalstrip}, we have
\begin{equation}\label{eq:Mbound}
 \left|\frac{a}{M_L(x+iy)}\right|
 \leq -\frac{a}{M_1(iy)}
 \leq -\frac{a}{M_1(i(1-\delta_a))}=1,
\end{equation}
for all $L\geq1$, $x\in\mathbb R$, and $|y|\leq 1-\delta_a$.
Consequently, the series
\[
 b_L(z)
 =\log\left(1+\frac{a}{M_L(z)}\right)
 =\sum_{n=1}^{\infty}
   \frac{(-1)^{n-1}a^n}{nM_L(z)^n}
\]
converges locally uniformly on the strip $|\Ima z|<1-\delta_a$ and
defines a single-valued holomorphic branch there. Since $M_L$ is
positive and even on the real axis, \eqref{eq:Mbound} shows that $b_L$
is real-valued and even there. On each closed substrip
$|\Ima z|\leq\eta<1-\delta_a$,
the strict form of \eqref{eq:Mbound} shows that $1+a/M_L(z)$ stays
uniformly away from zero. Together with \eqref{eq:MLinvsymbol}, this
gives \eqref{eq:blogbounds}. This completes the proof.
\end{proof}

\subsection{Analytic factorization}

We now construct the factor $w_L$ used in \eqref{eq:intro-multiplier}.
The convolution kernel in \eqref{eq:construction} is chosen to invert
the vertical averaging operator appearing in the shifted factorization
problem; its Fourier transform makes this mechanism explicit in
\eqref{eq:Ktransform}. Besides the boundary identities needed for
unitary conjugation, the construction gives estimates uniform in both
the spectral variable and $L$, which will be used in
Section~\ref{sec:multipliers}.

\begin{lemma}\label{lem:factor}
For each $L\geq1$, there exists a nonvanishing holomorphic function
$w_L$ on the strip $|\Ima z|<3/2-\delta_a$ such that 
\begin{align}
 |w_L(\xi)|&=1,\label{eq:unimodular}\\
 \frac{w_L(\xi+i/2)}{w_L(\xi-i/2)}
 &=\frac{M_L(\xi)+a}{M_L(\xi)},\label{eq:factor-ratio}\\
 |w_L(\xi-i/2)|^2&=\frac{M_L(\xi)}{M_L(\xi)+a}
 \label{eq:factoridentity}
\end{align}
for every $\xi\in\R$.
Moreover, for every $0\leq b<3/2-\delta_a$ and all nonnegative
integers $j,k$, there exists a constant $C_{a,b,j,k}>0$ such that
\begin{equation}\label{eq:wbound}
 |\partial_z^j\partial_L^k w_L(x+iy)|
 +|\partial_z^j\partial_L^k w_L^{-1}(x+iy)|
 \leq C_{a,b,j,k}L^{-k}(L+|x|)^{-j}
\end{equation}
for all $L\geq1$, $z=x+iy$, $x\in\mathbb R$, and $|y|\leq b$.
\end{lemma}

\begin{proof}
Let
\begin{equation}\label{eq:construction}
\begin{split}
 K(x)&=\frac{\pi}{2}\sech^2(\pi x),\qquad
 J_L=K*b_L,\\
 F_L(z)&=-i\int_0^z J_L(\zeta)\,d\zeta,\qquad
 w_L(z)=e^{F_L(z)}.
\end{split}
\end{equation}
On the real axis, \eqref{eq:construction} agrees with the definition
of $w_L$ in \eqref{WL}. We first establish the analytic
continuation and derivative bounds for $J_L$, and then verify
\eqref{eq:factor-ratio} and \eqref{eq:factoridentity}.

By analytic continuation of the convolution, $J_L$ is holomorphic
on $|\Ima z|<3/2-\delta_a$. More precisely,
\begin{equation}\label{eq:Jcontour}
 J_L(z)=\int_{\R}K(z-v-i\tau)b_L(v+i\tau)\,dv,
 \qquad |\tau|<1-\delta_a,\quad |\Ima z-\tau|<\frac12.
\end{equation}
The integral is independent of $\tau$, since $K$ has
exponential decay on closed substrips and no pole is crossed when
the contour is shifted. 

Fix $b<3/2-\delta_a$. For $|\Ima z|\leq b$, the contour in
\eqref{eq:Jcontour} can be chosen uniformly away from the boundary
of the admissible region. It follows from integration by parts and \eqref{eq:blogbounds} that 
\begin{equation}\label{eq:Jbound}
 |\partial_z^j\partial_L^kJ_L(x+iy)|
 \leq C_{a,b,j,k}(L+|x|)^{-1-j-k},
 \qquad |y|\leq b.
\end{equation}
Here we use
\[
 L+|x|\leq (L+|v|)(1+|x-v|)
\]
to control the convolution uniformly in $L$ and $x$.

Since $J_L$ is real on the real axis, $F_L$ is purely imaginary
there. Applying the Schwarz reflection symmetry of $J_L$, we have
\begin{equation}\label{eq:Freflection}
 F_L(\overline z)=-\overline{F_L(z)}.
\end{equation}
Integrating \eqref{eq:Jbound} along vertical segments, and along the
real axis for derivatives with respect to $L$, we obtain
\begin{align}
 |\Rea F_L(x+iy)|
 &\leq C_{a,b}(L+|x|)^{-1},
 \label{eq:Freal}\\
 |\partial_L^kF_L(x+iy)|
 &\leq C_{a,b,k}L^{-k},\qquad k\geq1,
 \label{eq:Fangular}\\
 |\partial_z^j\partial_L^kF_L(x+iy)|
 &\leq C_{a,b,j,k}(L+|x|)^{-j-k},
 \qquad j\geq1.
 \label{eq:Fderiv}
\end{align}
The chain rule applied to $w_L=e^{F_L}$ and
$w_L^{-1}=e^{-F_L}$ now yields \eqref{eq:wbound}.
Since $F_L$ is purely imaginary on $\R$ and $F_L(0)=0$,
we also obtain $w_L(0)=1$ and \eqref{eq:unimodular}.

It remains to establish the ratio identity. In the following calculation
we use the nonunitary Fourier transform $\widehat h(t)=\int e^{-itx}h(x)\,dx$.
The substitution $u=e^{2\pi x}$ and \cite{DLMF} give
\begin{equation}\label{eq:Ktransform}
 \widehat K(t)=\frac{t}{2\sinh(t/2)},\qquad
 \widehat K(0)=1.
\end{equation}

Thus $\widehat K(t)$ is the reciprocal of the multiplier
$\int_{-1/2}^{1/2}e^{-ty}\,dy$, with the value at $t=0$ understood
by continuity. This is the choice that turns vertical averaging of
$J_L$ into $b_L$.

We see that $J_L$ is holomorphic on a strip of width greater than
$1/2$ and satisfies \eqref{eq:Jbound}, which implies
\[
 \widehat{J_L(\cdot+iy)}(t)
 =e^{-ty}\widehat K(t)\widehat b_L(t),
 \qquad |y|\leq\frac12,
\]
in $L^2(\R)$. Consequently, we have
\[
 \widehat{\left(\int_{-1/2}^{1/2}J_L(\cdot+iy)\,dy\right)}(t)
 =\widehat K(t)
   \left(\int_{-1/2}^{1/2}e^{-ty}\,dy\right)\widehat b_L(t)
 =\widehat b_L(t),
\]
which yields
\[
 \int_{-1/2}^{1/2}J_L(x+iy)\,dy=b_L(x).
\]
It follows from $F_L'=-iJ_L$ that
\begin{equation}\label{eq:Fjump}
 F_L(x+i/2)-F_L(x-i/2)=b_L(x),
\end{equation}
which implies
\[
 \frac{w_L(x+i/2)}{w_L(x-i/2)}
 =e^{b_L(x)}
 =\frac{M_L(x)+a}{M_L(x)}
\]
by \eqref{eq:Fjump}.
Together with \eqref{eq:Freflection}, this means
\[
 |w_L(x-i/2)|^2
 =e^{2\Rea F_L(x-i/2)}
 =\frac{M_L(x)}{M_L(x)+a}.
\]
This completes the proof.
\end{proof}

\subsection{The unitary intertwiner}

We lift the scalar factors from Lemma~\ref{lem:factor} to the spatial
operator $W_a$ defined in \eqref{eq:W}. The shifted identity
\eqref{eq:factoridentity} is tailored to the weight $e^{-y/2}$ in the
form representation \eqref{eq:mellinform}; it therefore gives exact
conjugation, rather than merely comparison of the two forms.

\begin{proposition}\label{prop:intertwiner}
Let $\cW_a$ be defined by \eqref{eq:intro-cylinder-W}, and let
$W_a=\cT^{-1}\cW_a\cT$ as in \eqref{eq:W}. Then $W_a$ is unitary
on $L^2(\R^3)$ and maps $\Dom(q_0)$ onto $\Dom(q_a)$. Moreover,
\begin{equation}\label{eq:structure}
 H_a=W_a|D|W_a^*,\qquad
 m(H_a)=W_am(|D|)W_a^*.
\end{equation}
The second identity holds for Borel functions $m$, with both
sides defined on their respective spectral domains.

\end{proposition}

\begin{proof}
By \eqref{eq:unimodular}, each $w_{\ell+1}(D_y)$ is unitary on
$L^2(\R)$; hence \eqref{eq:intro-cylinder-W} and \eqref{eq:W} imply
that $W_a$ is unitary on $L^2(\R^3)$.

If $v\in C_c^\infty(\R)$ and $\beta\in\R$ satisfies
$|\beta|<3/2-\delta_a$, then
\begin{equation}\label{eq:shift}
 e^{-\beta y}w_L(D_y)^{\pm1}v
 =w_L(D_y-i\beta)^{\pm1}(e^{-\beta y}v).
\end{equation}
Indeed, we note that
\[
 \cF_y(e^{-\beta y}v)(\xi)
 =\widehat v(\xi-i\beta).
\]
Thus \eqref{eq:shift} follows by shifting the Fourier contour.
The rapid decay of $\widehat v$ on closed horizontal strips and
the bounds in \eqref{eq:wbound} justify the contour shift, since
the integrals over the vertical sides vanish as their real parts
tend to $\pm\infty$.

Since $3/2-\delta_a>1/2$, we may take $\beta=1/2$ in
\eqref{eq:shift}. This is precisely the shift required by the weight
in \eqref{eq:mellinform}. Let $\mathcal C$ denote the common form core
from Lemma~\ref{lem:mellinform}. For $f\in\mathcal C$, write
\[
 \cT f=\sum_{\ell,m}u_{\ell m}(y)Y_{\ell m}(\omega).
\]
According to Lemma~\ref{lem:mellinform},
together with \eqref{eq:factoridentity} and \eqref{eq:shift},
we obtain
\begin{align*}
 q_a[W_af]
 &=\sum_{\ell,m}\int_{\R}
 (M_{\ell+1}(\xi)+a)
 |w_{\ell+1}(\xi-i/2)|^2
 \left|\cF_y(e^{-y/2}u_{\ell m})(\xi)\right|^2
 \,\dd\xi\\
 &=\sum_{\ell,m}\int_{\R}
 M_{\ell+1}(\xi)
 \left|\cF_y(e^{-y/2}u_{\ell m})(\xi)\right|^2
 \,\dd\xi\\
 &=q_0[f].
\end{align*}
Similarly, applying the same argument to the reciprocal
multipliers gives
\begin{align*}
 q_0[W_a^*f]
 &=\sum_{\ell,m}\int_{\R}
 M_{\ell+1}(\xi)
 |w_{\ell+1}(\xi-i/2)|^{-2}
 \left|\cF_y(e^{-y/2}u_{\ell m})(\xi)\right|^2
 \,\dd\xi\\
 &=q_a[f].
\end{align*}
The maximal-domain characterization in Lemma~\ref{lem:mellinform}
ensures that these identities hold on the indicated form cores
and that the images belong to the corresponding form domains.

Since $\mathcal C$ is a form core for both $q_0$ and $q_a$,
the preceding identities extend by approximation in the
respective form norms to
\[
 W_a:\Dom(q_0)\longrightarrow\Dom(q_a),
 \qquad
 W_a^*:\Dom(q_a)\longrightarrow\Dom(q_0).
\]
Since $W_a$ and $W_a^*$ are unitary operators, $W_a$ maps
$\Dom(q_0)$ onto $\Dom(q_a)$.

It follows from polarization and the representation theorem for closed quadratic forms that $H_a=W_a|D|W_a^*$. The Borel functional calculus
then yields
\[
 m(H_a)=W_am(|D|)W_a^*
\]
on the corresponding spectral domains.
\end{proof}

\section{\texorpdfstring{$L^p$}{Lp} bounds for the intertwining operators}
\label{sec:multipliers}

This section completes the proof of the first main result,
Theorem~\ref{thm:intertwining}. We combine the mixed estimates
\eqref{eq:wbound} for the analytic factors with a joint multiplier
theorem on the logarithmic cylinder. The resulting argument preserves
uniformity in the angular momentum and identifies the full interval
$I_a$ in \eqref{eq:interval}.

\subsection{Joint multipliers on the logarithmic cylinder}

The next lemma turns estimates that are simultaneous in the radial
frequency $\xi$ and the angular parameter $L$ into an
$L^p(\R\times\Sph)$ multiplier theorem. It is the device that allows
the angular channels to be recombined with constants uniform in
$\ell$.

\begin{lemma}[Joint multiplier theorem on the cylinder]
\label{lem:product}
Let $m\in C^\infty(\R\times[1,\infty))$ satisfy
\begin{equation}\label{eq:productbounds}
 \sup_{\xi\in\R,\,L\geq1}
 |\xi|^jL^k
 |\partial_\xi^j\partial_L^k m(\xi,L)|<\infty
\end{equation}
for all nonnegative integers $j,k$. Define the operator
$\mathcal M_m$ on functions with finitely many spherical harmonic
components by
\[
 \cF_y\bigl[(\mathcal M_m f)_{\ell n}\bigr](\xi)
 =m(\xi,\ell+1)\cF_y f_{\ell n}(\xi),
\]
where $f_{\ell n}$ is the coefficient in the spherical harmonic
expansion \eqref{eq:intro-expansion}.

Then $\mathcal M_m$ extends to a bounded operator on
$L^p(\R\times\Sph)$ for every $1<p<\infty$. 
\end{lemma}

\begin{proof}
Consider the strongly commuting operators
\[
 B_1=|D_y|\otimes I,\qquad
 B_2=I\otimes\Lambda,
\]
where
\[
 \Lambda=\sqrt{-\Delta_{\Sph}+\frac14}+\frac12.
\]

Both $B_1$ and $B_2$ generate symmetric contraction semigroups.
For $B_1$, this is the Poisson semigroup on $\R$. For $B_2$, let
\[
 L_0=-\Delta_{\Sph}+\frac14.
\]
Then
\[
 e^{-tB_2}=e^{-t/2}e^{-t\sqrt{L_0}},
\]
where the second factor is obtained by subordination from the
heat semigroup generated by $L_0$. Since the Fourier transform of an $L^2(\R)$ function cannot be supported at the single point $\xi=0$, we have $E_{B_1}(\{0\})=0$. It follows from $\Lambda\geq1$ that $E_{B_2}(\{0\})=0$.

We next verify that both operators admit Marcinkiewicz functional
calculi of finite order. For $B_1$, the required annular $L^2$
derivative estimates follow from the one-dimensional Fourier
multiplier theorem \cite{Stein}. For $B_2$, we apply the spectral multiplier theorem of Duong, Sikora, and Yan \cite{DuongSikoraYan}  to $L_0$.
The round sphere is a doubling space of dimension two, and the
heat kernel of $L_0$ satisfies Gaussian bounds. Hence a uniform
localized Sobolev condition of order three in $L^\infty$ is
sufficient for $L^p$ boundedness for every $1<p<\infty$.

To verify this condition, let $v$ be a smooth scalar multiplier
of $\Lambda$ and define
\[
 F(\lambda)=v\left(\sqrt{\lambda}+\frac12\right),
 \qquad \lambda\geq0.
\]
Then $v(\Lambda)=F(L_0)$. The chain rule and one-dimensional
Sobolev embedding show that annular $L^2$ bounds on five
derivatives of $v$ imply the required localized Sobolev bounds
for $F$. The low-frequency value $F(0)=v(1/2)$ is controlled
by the same estimates after a smooth extension of $v$ below
$L=1$. Thus $B_2$ also
admits a Marcinkiewicz functional calculus of finite order.
These one-variable estimates remain valid on the product space
$L^p(\R\times\Sph)$ by Fubini's theorem.

By Wr\'obel's joint calculus theorem \cite{Wrobel}, individual orders $\rho_1,\rho_2$ give joint orders strictly larger than $(\rho_1+1,\rho_2+1)$. The pointwise bounds in \eqref{eq:productbounds} imply all the
mixed annular $L^2$ derivative estimates required by this
theorem. 

Finally, we handle the sign of the radial Fourier variable.
Let $\Pi_\pm$ denote the projections onto
positive and negative Fourier frequencies in $y$, which are bounded on $L^p(\R\times\Sph)$ for
$1<p<\infty$ by the one-dimensional multiplier theorem. Define $m_\pm(\lambda,L)=m(\pm\lambda,L)$ for $\lambda>0$, which 
satisfy the same mixed derivative bounds as $m$.
Because the spectrum of $B_2$ is contained in $[1,\infty)$, we may
extend $m_\pm$ smoothly to $(0,\infty)^2$ in the $L$ variable without
changing their values on the joint spectrum. A fixed extension below
$L=1$ preserves the required Marcinkiewicz bounds. No value at
$\lambda=0$ is needed because $E_{B_1}(\{0\})=0$.
We have
\[
 \mathcal M_m
 =\Pi_+\mathcal M_{m_+}(B_1,B_2)
  +\Pi_-\mathcal M_{m_-}(B_1,B_2).
\]
According to the boundedness of $\Pi_\pm$ and Wr\'obel's joint calculus theorem \cite{Wrobel}, $\mathcal M_m$ extends to a bounded operator on
$L^p(\R\times\Sph)$ for every $1<p<\infty$.
\end{proof}

\subsection{\texorpdfstring{$L^p$}{Lp} bounds for the intertwiner}

We apply Lemma~\ref{lem:product} to the factors from
Lemma~\ref{lem:factor}. The notable point is that the $L^p$ polar
isometry produces an imaginary shift of the Mellin variable. Hence the
holomorphic strip in Lemma~\ref{lem:factor} determines exactly the
interval $I_a$ in \eqref{eq:interval}.

\begin{proof}[Proof of Theorem~\ref{thm:intertwining}]
The construction, unitarity, and conjugation identity were proved in
Proposition~\ref{prop:intertwiner}. It remains to prove the $L^p$ bounds.
Fix $p\in I_a$, and let
\[
 \beta_p=\frac32-\frac3p.
\]
Define
\begin{equation}\label{eq:Tp}
 (\cT_p f)(y,\omega)
 =e^{3y/p}f(e^y\omega)
 =e^{-\beta_p y}\cT f(y,\omega).
\end{equation}
This is an isometry from $L^p(\R^3)$ onto $L^p(\R\times\Sph)$.
On the common core from Lemma~\ref{lem:mellinform},
\eqref{eq:shift} gives
\begin{equation}\label{eq:Lpshift}
 \cT_p W_a^{\pm1}f
 =w_\Lambda(D_y-i\beta_p)^{\pm1}\cT_p f.
\end{equation}
Since $p\in I_a$, we have $|\beta_p|<3/2-\delta_a$. Choose
$b\in[0,3/2-\delta_a)$ such that $|\beta_p|\leq b$. By
\eqref{eq:wbound}, for all nonnegative integers $j,k$,
\[
 |\xi|^jL^k
 \left|\partial_\xi^j\partial_L^k
 w_L^{\pm1}(\xi-i\beta_p)\right|
 \leq
 C_{a,b,j,k}
 \left(\frac{|\xi|}{L+|\xi|}\right)^j
 \leq C_{a,b,j,k},
 \qquad \xi\in\R,\quad L\geq1.
\]
Lemma~\ref{lem:product} therefore shows that
$w_\Lambda(D_y-i\beta_p)^{\pm1}$ is bounded on
$L^p(\R\times\Sph)$. In view of \eqref{eq:Lpshift}, $W_a$ and
$W_a^*$ are bounded inverse operators on $L^p(\R^3)$.
\end{proof}

\section{Strichartz estimates and sharpness}
\label{sec:strichartz}

We turn to the second main result, Theorem~\ref{thm:main}. Its proof
requires two additional ingredients. We first compare the free and
perturbed Sobolev norms in the range dictated by the analytic strip
from Lemma~\ref{lem:factor}. We then transfer the free half-wave
estimates through the intertwining identity \eqref{eq:structure}.
The section concludes by proving that the spatial restriction in the
attractive case is necessary.

\subsection{A one-sided Sobolev comparison}

We first record a comparison valid throughout the
coupling-independent range $0\leq\sigma<3/2$. For
$1<\sigma<3/2$, the proof bootstraps the known lower-order estimate
using a Stein--Weiss bound for multiplication by the Coulomb
potential; it does not commute fractional powers through $H_a$.

\begin{proposition}[One-sided Sobolev comparison]
\label{prop:upper-sobolev}
Let $a>-2/\pi$. For every $0\leq\sigma<3/2$, we have
\begin{equation}
 H^\sigma(\R^3)\subset\Dom(H_a^\sigma),\qquad
 \norm{H_a^\sigma f}_{L^2}
 \leq C_{a,\sigma}\norm{|D|^\sigma f}_{L^2},
 \quad f\in H^\sigma(\R^3).
 \label{eq:upper-sobolev}
\end{equation}
Consequently, the operator $H_a^\sigma$, initially defined on
$C_c^\infty(\R^3\setminus\{0\})$, extends uniquely to a bounded
operator from $\dot H^\sigma(\R^3)$ to $L^2(\R^3)$.
\end{proposition}

\begin{proof}
The case $\sigma=0$ is immediate. For $0<\sigma\leq1$, 
both the norm estimate and the operator-domain inclusion are
established in \cite{FMS}. Indeed, the estimate follows by
applying that result with $d=3$, $\alpha=1$, and power parameter
$2\sigma$. 

For $1<\sigma<3/2$, let $\theta=\sigma-1\in(0,1/2)$. 
It follows from the Stein--Weiss weighted fractional integral
inequality \cite{SteinWeiss} that 
\begin{equation}\label{eq:coulomb-multiplication}
 \norm{|D|^\theta(|x|^{-1}f)}_{L^2}
 \leq C_\theta\norm{|D|^{1+\theta}f}_{L^2}.
\end{equation}
Indeed, the homogeneous Fourier transform formula \eqref{eq:homogeneous} gives
$\widehat{|x|^{-1}}(\xi)=c|\xi|^{-2}$ for a positive constant
$c$. Consequently, after setting
$h(\eta)=|\eta|^{1+\theta}\widehat f(\eta)$, the relevant
Fourier-space operator, up to a constant factor, is
\[
 h\longmapsto
 |\xi|^\theta
 \int_{\R^3}
 \frac{|\eta|^{-1-\theta}h(\eta)}
      {|\xi-\eta|^2}\,\dd\eta.
\]
The Stein--Weiss parameters are $1+\theta$ on the input and
$-\theta$ on the output, and their sum is $1$, as required
for an order-one fractional integral. The strict admissibility
conditions hold because $1+\theta<3/2$ and $-\theta<3/2$.
This proves \eqref{eq:coulomb-multiplication}.

Let $f\in C_c^\infty(\R^3\setminus\{0\})$. 
Estimate \eqref{eq:coulomb-multiplication} implies that
$H_af\in H^\theta(\R^3)$. 
It follows from $H^\theta(\R^3)\subset\Dom(H_a^\theta)$ that 
$H_af\in\Dom(H_a^\theta)$.
Here we use the already established case $0<\theta<1/2$ of
\eqref{eq:upper-sobolev}. Since
$f\in C_c^\infty(\R^3\setminus\{0\})\subset\Dom(H_a)$, the spectral
theorem then yields
$f\in\Dom(H_a^{1+\theta})$ and
\begin{align}
 \norm{H_a^{1+\theta}f}_{L^2}
 &=\norm{H_a^\theta H_af}_{L^2} \notag\\
 &\leq C_{a,\theta}\norm{|D|^\theta H_af}_{L^2} \notag\\
 &\leq C_{a,\theta}
 \left(
   \norm{|D|^{1+\theta}f}_{L^2}
   +|a|\norm{|D|^\theta(|x|^{-1}f)}_{L^2}
 \right) \notag\\
 &\leq C_{a,\theta}\norm{|D|^{1+\theta}f}_{L^2}.
 \label{eq:sobolev-bootstrap}
\end{align}
Thus the desired estimate follows for
$1<\sigma<3/2$, initially on
$C_c^\infty(\R^3\setminus\{0\})$. In particular, the argument
uses the known comparison at the lower order $\theta$ and
\eqref{eq:coulomb-multiplication}, without commuting fractional
powers through $H_a$.

For every $0\leq\sigma<3/2$, the space
$C_c^\infty(\R^3\setminus\{0\})$ is dense in both
$H^\sigma(\R^3)$ and $\dot H^\sigma(\R^3)$.
To see this, choose
$\chi\in C_c^\infty(\R^3)$ with $\chi=1$ near the origin.
For $f\in C_c^\infty(\R^3)$, scaling gives
\begin{equation}
\begin{aligned}
 \norm{|D|^\sigma[\chi(x/\varepsilon)f(x)]}_{L^2} &=
 \varepsilon^{3/2-\sigma}
 \norm{|D|^\sigma[\chi(\cdot)f(\varepsilon\cdot)]}_{L^2} \\
 &\leq C_f\varepsilon^{3/2-\sigma}
 \longrightarrow0
\end{aligned}
\label{eq:punctured-density}
\end{equation}
as $\varepsilon\to0$. The final inequality follows by
interpolation between $L^2$ and two derivatives: the rescaled
functions have uniformly bounded derivatives through order two
and are supported in a fixed compact set. The removed portions
also tend to zero in $L^2$. This establishes density in the
inhomogeneous Sobolev norm as well.

Applying the estimate already proved on the common dense core,
we may pass to the limit in the graph norm of $H_a^\sigma$.
Since $H_a^\sigma$ is closed, this proves
\eqref{eq:upper-sobolev} on $H^\sigma(\R^3)$. Approximation in
the homogeneous Sobolev norm likewise yields a unique bounded
extension from $\dot H^\sigma(\R^3)$ to $L^2(\R^3)$. 
\end{proof}

\subsection{Homogeneous Sobolev norms on the cylinder}

We next express the Euclidean homogeneous Sobolev norm in the same
logarithmic polar coordinates used to construct $W_a$. The exact
identity below, rather than only a norm equivalence, makes the
contour-shift argument compatible with the joint multipliers
$w_\Lambda(D_y)$.

Let $\cT$ be as in \eqref{eq:intro-polar}, let $D_y=-i\partial_y$,
and recall that $\Lambda Y_{\ell m}=(\ell+1)Y_{\ell m}$. For
$0<\sigma<3/2$, define
\begin{equation}\label{eq:dsigma}
 d_{\sigma,\ell}(\xi)
 :=2^{2\sigma}
 \left|
 \frac{
 \Gamma\bigl((\ell+3/2+\sigma+i\xi)/2\bigr)
 }{
 \Gamma\bigl((\ell+3/2-\sigma+i\xi)/2\bigr)
 }
 \right|^2,
 \qquad d_{0,\ell}(\xi):=1,
\end{equation}
and let $D_\sigma(D_y,\Lambda)$ denote the corresponding joint
Fourier--spectral multiplier. Applying
\eqref{eq:gamma-ratio-asymptotic} with
\[
 w=\frac{\ell+3/2-\sigma+i\xi}{2},
 \qquad (\alpha,\beta)=(\sigma,0),
\]
gives
\begin{equation}\label{eq:dsigma-asymptotics}
 d_{\sigma,\ell}(\xi)
 \asymp_\sigma (\ell+1+|\xi|)^{2\sigma},
 \qquad \ell\geq0,\quad \xi\in\R.
\end{equation}
Indeed, the real part of the denominator argument satisfies
\[
 \frac{\ell+3/2-\sigma}{2}
 \geq \frac{3/2-\sigma}{2}>0.
\]
Thus, for each fixed $\sigma\in(0,3/2)$, both Gamma-function
arguments remain in a right half-plane, and the uniform
Gamma-quotient asymptotics apply. The estimate in the remaining
bounded region follows from continuity and the fact that the
Gamma function has no zeros.

\begin{lemma}\label{lem:sobolev-cylinder}
Let $0\leq\sigma<3/2$. Then, for every
$f\in H^\sigma(\R^3)$,
\begin{equation}
 \norm{|D|^\sigma f}_{L^2}
 =
 \norm{
 D_\sigma(D_y,\Lambda)^{1/2}
 e^{-\sigma y}\cT f
 }_{L^2}.
 \label{eq:sobolev-cylinder}
\end{equation}
For $0<\sigma<3/2$, the image of $H^\sigma(\R^3)$ under
$\cT$ is precisely
\begin{equation}
 \mathscr D_\sigma
 =
 \left\{
 F\in L^2(\R\times\Sph):
 e^{-\sigma y}F
 \in\Dom\bigl(D_\sigma(D_y,\Lambda)^{1/2}\bigr)
 \right\}.
 \label{eq:sobolev-cylinder-domain}
\end{equation}
\end{lemma}

\begin{proof}
We first establish \eqref{eq:sobolev-cylinder} on a dense class
of smooth functions supported in a compact annulus and having
finitely many spherical harmonic components.

The Fourier transform of a homogeneous spherical harmonic
\cite{SteinWeissBook} yields
\begin{equation}\label{eq:homogeneous-action-sobolev}
 |D|^{2\sigma}
 \bigl(r^{-3/2+\sigma+i\xi}Y_{\ell m}(\omega)\bigr)
 =
 d_{\sigma,\ell}(\xi)
 r^{-3/2-\sigma+i\xi}Y_{\ell m}(\omega).
\end{equation}
Here $d_{\sigma,\ell}$ is defined in
\eqref{eq:dsigma}. The restriction
$\sigma<3/2$ ensures that the homogeneous distributions
appearing in \eqref{eq:homogeneous-action-sobolev} are
represented by locally integrable functions.

Let
\[
 F=\cT f,
 \qquad
 G=e^{-\sigma y}F.
\]
Since
\[
 f(r\omega)=r^{-3/2+\sigma}G(\log r,\omega),
\]
Mellin inversion and \eqref{eq:homogeneous-action-sobolev}
give
\[
 (|D|^{2\sigma}f)(r\omega)
 =
 r^{-3/2-\sigma}
 \bigl(D_\sigma(D_y,\Lambda)G\bigr)(\log r,\omega)
\]
in the distributional sense. The uniform bounds justify the
superposition of the homogeneous modes for the stated class
of functions.

Pairing this identity with $f$ and using polar coordinates,
we obtain
\begin{align*}
 \norm{|D|^\sigma f}_{L^2}^2
 &=\ip{|D|^{2\sigma}f}{f}_{L^2(\R^3)}\\
 &=\ip{D_\sigma(D_y,\Lambda)G}{G}_{L^2(\R\times\Sph)}\\
 &=\norm{D_\sigma(D_y,\Lambda)^{1/2}G}_{L^2}^2.
\end{align*}
This proves \eqref{eq:sobolev-cylinder} on the stated
dense class.

For $0<\sigma<3/2$, we note that
\begin{equation}\label{eq:cylinder-form-equivalence}
 \norm{D_\sigma(D_y,\Lambda)^{1/2}G}_{L^2}
 \asymp_\sigma
 \norm{\Lambda^\sigma G}_{L^2}
 +\norm{|D_y|^\sigma G}_{L^2}.
\end{equation}
Consequently, the domain of
$D_\sigma(D_y,\Lambda)^{1/2}$ is
$H^\sigma(\R\times\Sph)$, with an equivalent norm.
In particular,
\[
 \norm{G}_{L^2}
 \leq C_\sigma
 \norm{D_\sigma(D_y,\Lambda)^{1/2}G}_{L^2}.
\]

Consider the quadratic form
\[
 q_\sigma[F]
 =
 \norm{D_\sigma(D_y,\Lambda)^{1/2}
 e^{-\sigma y}F}_{L^2}^2
\]
on $\mathscr D_\sigma$. We claim that this form is closed.
Suppose that $F_n\to F$ in $L^2(\R\times\Sph)$ and
\[
 D_\sigma(D_y,\Lambda)^{1/2}e^{-\sigma y}F_n
 \longrightarrow H
 \quad\text{in }L^2(\R\times\Sph).
\]
By \eqref{eq:cylinder-form-equivalence},
$e^{-\sigma y}F_n$ is a Cauchy sequence in $L^2$.
Moreover, multiplication by $e^{-\sigma y}$ is continuous
on every compact interval in $y$. Hence its limit is
$e^{-\sigma y}F$. The closedness of
$D_\sigma(D_y,\Lambda)^{1/2}$ now gives
\[
 e^{-\sigma y}F
 \in\Dom\bigl(D_\sigma(D_y,\Lambda)^{1/2}\bigr),
 \qquad
 D_\sigma(D_y,\Lambda)^{1/2}e^{-\sigma y}F=H.
\]
Thus $q_\sigma$ is closed.

Next, we prove that smooth functions compactly supported in $y$ and having finitely many spherical harmonic components form a core for the corresponding quadratic form. Let
$F\in\mathscr D_\sigma$ and set $G=e^{-\sigma y}F$.
By \eqref{eq:cylinder-form-equivalence},
$G\in H^\sigma(\R\times\Sph)$. Choose smooth cutoffs
$\chi_R(y)=\chi(y/R)$, where
$\chi\in C_c^\infty(\R)$ and $\chi=1$ near the origin.
Then
\[
 \chi_R F\longrightarrow F
 \quad\text{in }L^2,
 \qquad
 \chi_R G\longrightarrow G
 \quad\text{in }H^\sigma
\]
as $R\to\infty$. The second convergence follows from the
standard density property of smooth cutoffs in fractional
Sobolev spaces. For $0<\sigma<1$, it follows directly from
the difference-quotient characterization of $H^\sigma$;
for $1\leq\sigma<3/2$, it follows by applying the same
argument to first derivatives and using the product rule.

For each fixed $R$, we may approximate $\chi_R G$ in
$H^\sigma(\R\times\Sph)$ by smooth functions $G_n$ compactly
supported in $y$ and having finitely many spherical harmonic
components. The approximants can be chosen with support in a
fixed compact interval. Setting
\[
 F_n=e^{\sigma y}G_n,
\]
we obtain smooth functions of the same type and
\[
 F_n\to\chi_R F\quad\text{in }L^2,
 \qquad
 e^{-\sigma y}F_n\to\chi_R G
 \quad\text{in }H^\sigma.
\]
A diagonal argument proves the asserted core property.

The preceding calculation therefore identifies the two
quadratic forms on a common core. Both forms are closed:
the Euclidean form has domain $H^\sigma(\R^3)$, while the
cylinder form has domain $\mathscr D_\sigma$. Since
$\cT$ is unitary on $L^2$ and maps the common core onto the
corresponding cylinder core, the uniqueness of the closures
of these forms yields
\[
 \cT\bigl(H^\sigma(\R^3)\bigr)=\mathscr D_\sigma
\]
and extends \eqref{eq:sobolev-cylinder} to every
$f\in H^\sigma(\R^3)$.

Finally, when $\sigma=0$, we have $D_0=I$ and
$\cT$ is unitary. Thus \eqref{eq:sobolev-cylinder}
holds trivially.
\end{proof}

\subsection{Sobolev isomorphisms}

We now combine the cylinder identity \eqref{eq:sobolev-cylinder} with
the strip estimates \eqref{eq:wbound}. The same analytic strip that
determines the $L^p$ interval in \eqref{eq:interval} also determines
the Sobolev range below. Thus the two mapping ranges have a common
origin in the factorization from Lemma~\ref{lem:factor}.

\begin{lemma}
\label{lem:sobolev-isomorphism}
Let $a>-2/\pi$ and $0\leq\sigma<3/2-\delta_a$. Then $W_a$ and
$W_a^*$ extend to bounded operators on $H^\sigma(\R^3)$ and
to mutually inverse bounded isomorphisms on $\dot H^\sigma(\R^3)$.
Moreover,
\begin{equation}
 \Dom(H_a^\sigma)=H^\sigma(\R^3),\qquad
 \norm{H_a^\sigma f}_{L^2}
 \asymp_{a,\sigma}\norm{|D|^\sigma f}_{L^2}.
 \label{eq:two-sided-sobolev}
\end{equation}
\end{lemma}

\begin{proof}
The assertion is immediate for $\sigma=0$, since $W_a$ is unitary
on $L^2(\R^3)$. We first establish the Sobolev isomorphism for
$0<\sigma<3/2-\delta_a$.

By Lemma~\ref{lem:sobolev-cylinder}, the homogeneous Sobolev norm
admits the cylinder representation
\[
 \norm{|D|^\sigma f}_{L^2}
 =
 \norm{D_\sigma(D_y,\Lambda)^{1/2}
 e^{-\sigma y}\cT f}_{L^2}.
\]
On the common core of smooth functions compactly supported in $y$
and having finitely many spherical harmonic components, the
intertwining identity and the contour-shift relation give
\begin{equation}
 e^{-\sigma y}w_\Lambda(D_y)^{\pm1}F
 =
 w_\Lambda(D_y-i\sigma)^{\pm1}e^{-\sigma y}F.
 \label{eq:sobolev-shift}
\end{equation}
Lemma~\ref{lem:factor} implies that the
shifted multipliers
$w_\Lambda(D_y-i\sigma)^{\pm1}$ are bounded on
$L^2(\R\times\Sph)$. They also commute with
$D_\sigma(D_y,\Lambda)^{1/2}$, as both are joint spectral
multipliers of $D_y$ and $\Lambda$.

It follows from \eqref{eq:sobolev-shift} that
\begin{equation}
 \norm{|D|^\sigma W_a^{\pm1}f}_{L^2}
 \leq M_{a,\sigma}\norm{|D|^\sigma f}_{L^2}
 \label{eq:W-sobolev}
\end{equation}
on the common core. By density, this estimate extends to
$\dot H^\sigma(\R^3)$. Since $W_a$ and $W_a^*$ are also unitary
on $L^2(\R^3)$, they extend to bounded operators on
$H^\sigma(\R^3)$. Their inverse identities persist under these
extensions, so they define mutually inverse isomorphisms on both
Sobolev spaces.

By the functional-calculus identity \eqref{eq:structure},
\[
 H_a^\sigma=W_a|D|^\sigma W_a^*.
\]
Consequently,
\[
 \Dom(H_a^\sigma)
 =W_a\Dom(|D|^\sigma)
 =W_aH^\sigma(\R^3)
 =H^\sigma(\R^3).
\]
Moreover,
\[
 \norm{H_a^\sigma f}_{L^2}
 =\norm{|D|^\sigma W_a^*f}_{L^2}.
\]
Applying \eqref{eq:W-sobolev} to $W_a^*$ and using the inverse
isomorphism property yields \eqref{eq:two-sided-sobolev}.
\end{proof}

\begin{remark}
The range $\sigma<3/2-\delta_a$ in
Lemma~\ref{lem:sobolev-isomorphism}
comes from the available holomorphic strip. The lemma makes no claim
about the endpoint. The one-sided estimate \eqref{eq:upper-sobolev}
holds throughout $0\leq\sigma<3/2$, including the range
$3/2-\delta_a\leq\sigma<3/2$ when $-2/\pi<a<0$.
\end{remark}

\subsection{Transfer of the free half-wave estimates}

We begin with the free estimate of Keel and Tao \cite{KT}.

\begin{proposition}[Free half-wave Strichartz estimates]
\label{prop:free-strichartz}
For every pair $(r,p)$ satisfying
\[
 2<r\leq\infty,\qquad
 2\leq p<\infty,\qquad
 \frac1r+\frac1p\leq\frac12,
\]
let
\[
 s=\frac32-\frac3p-\frac1r.
\]
Then
\begin{equation}\label{eq:free-strichartz}
 \norm{e^{-it|D|}g}_{L_t^rL_x^p}
 \leq C_{r,p}\norm{|D|^s g}_{L^2}
\end{equation}
for every $g\in\dot H^s$.
\end{proposition}

We now combine Proposition~\ref{prop:free-strichartz} with the
$L^p$ boundedness from Theorem~\ref{thm:intertwining} and the Sobolev
comparison \eqref{eq:two-sided-sobolev}.

\begin{proof}[Proof of Theorem~\ref{thm:main}]
Fix a pair $(r,p)$ satisfying \eqref{eq:admissible}
and let
\[
 s=\frac32-\frac1r-\frac3p.
\]
Suppose that $f\in\Dom(H_a^s)$ and $g=W_a^*f$. By the intertwining
identity \eqref{eq:structure},
\[
 e^{-itH_a}f=W_ae^{-it|D|}g,
 \qquad
 \norm{|D|^sg}_{L^2}=\norm{H_a^sf}_{L^2}.
\]
Since $p\in I_a$, Theorem~\ref{thm:intertwining} gives the boundedness
of $W_a$ on $L^p(\R^3)$. Applying this bound and then
\eqref{eq:free-strichartz}, we obtain
\begin{equation}\label{eq:spectral-strichartz}
 \norm{e^{-itH_a}f}_{L_t^rL_x^p}
 \leq C_{a,r,p}\norm{H_a^sf}_{L^2}.
\end{equation}
The admissibility conditions imply $0\leq s<3/2-\delta_a$:
indeed $1/r+3/p\leq3/2$, while $3/p>\delta_a$.
It follows from Lemma~\ref{lem:sobolev-isomorphism} that
\[
 \norm{H_a^sf}_{L^2}\leq C_{a,s}\norm{|D|^sf}_{L^2},
 \qquad f\in H^s(\R^3).
\]
Combining this estimate with \eqref{eq:spectral-strichartz}, we have
\begin{equation}\label{eq:ordinary-strichartz}
 \norm{e^{-itH_a}f}_{L_t^rL_x^p}
 \leq C_{a,r,p}\norm{|D|^sf}_{L^2},
 \qquad f\in H^s.
\end{equation}
Since $s<3/2$, smooth functions supported away from the origin are
dense in $\dot H^s(\R^3)$ by \eqref{eq:punctured-density}. Hence
\eqref{eq:ordinary-strichartz} defines a unique continuous extension
of the evolution to every $f\in\dot H^s(\R^3)$. On
$\dot H^s\cap L^2$, this extension agrees with the spectral-theoretic
evolution $e^{-itH_a}f$, which completes the proof.
\end{proof}

\subsection{Sharpness in the attractive case}

We finally prove that the restriction $p<3/\delta_a$ in
Theorem~\ref{thm:main} is necessary when $a<0$. The obstruction is the
ground-state singularity visible in the heat-kernel lower bound; this
links the dispersive range directly to the same exponent $\delta_a$
that controls the analytic strip in Lemma~\ref{lem:factor}.

\begin{proposition}\label{prop:attractive-obstruction}
Let $-2/\pi<a<0$. If $p\geq3/\delta_a$,
then, for every $1\leq r\leq\infty$ and every $s\geq0$, there is
no constant $C>0$ such that
\begin{equation}\label{eq:impossible-strichartz}
 \norm{e^{-itH_a}f}_{L_t^rL_x^p}
 \leq C\norm{|D|^sf}_{L^2}
\end{equation}
for all $f\in C_c^\infty(\R^3\setminus\{0\})$.
\end{proposition}

\begin{proof}
We first show that the heat semigroup fails to preserve $L^p$
under the stated condition on $p$. The two-sided heat-kernel
estimates of Bogdan, Grzywny, Jakubowski, and Pilarczyk
\cite{BGJP}, specialized to dimension $3$ and fractional order
$1$, imply that
\begin{equation}\label{eq:heat-lower-obstruction}
 e^{-H_a}(x,y)\geq c|x|^{-\delta_a},
 \qquad 0<|x|<\frac12,\quad 1<|y|<2.
\end{equation}

Choose a nonnegative, nonzero function
$f\in C_c^\infty(\R^3)$ supported in the annulus
$\{y:1<|y|<2\}$. Integrating \eqref{eq:heat-lower-obstruction}
against $f$ gives
\[
 e^{-H_a}f(x)\geq c_f|x|^{-\delta_a},
 \qquad 0<|x|<\frac12,
\]
for some constant $c_f>0$. If $p<\infty$, it follows that
\[
 \norm{e^{-H_a}f}_p^p
 \geq c_f^p\int_{|x|<1/2}|x|^{-p\delta_a}\,\dd x
 =\infty,
\]
because
\[
 \int_0^{1/2}r^{2-p\delta_a}\,\dd r=\infty
 \qquad\text{when }p\delta_a\geq3.
\]
At the endpoint $p\delta_a=3$, the divergence is logarithmic.
If $p=\infty$, the same pointwise lower bound shows that
$e^{-H_a}f$ is essentially unbounded near the origin.
Consequently,
\begin{equation}\label{eq:heat-not-Lp}
 e^{-H_a}f\notin L^p(\R^3).
\end{equation}

Suppose, to the contrary, that \eqref{eq:impossible-strichartz}
holds for some $r\in[1,\infty]$ and $s\geq0$. Choose
$m\in\mathcal S(\R)$ such that
\[
 m(\lambda)=e^{-\lambda},\qquad \lambda\geq0,
\]
and let $g\in\mathcal S(\R)$ be its inverse Fourier transform,
normalized so that
\[
 m(\lambda)=\int_{\R}g(t)e^{-it\lambda}\,\dd t.
\]
Since $H_a$ is nonnegative, the spectral theorem yields
\begin{equation}\label{eq:heat-fourier-representation}
 e^{-H_a}f
 =\int_{\R}g(t)e^{-itH_a}f\,\dd t
\end{equation}
as an identity in $L^2(\R^3)$. Under the assumed spacetime estimate,
the same integral converges in $L^p$ when $p<\infty$; equivalently, it
may first be interpreted weakly against $L^{p'}\cap L^2$ functions.

For finite $p$, the assumed estimate and H\"older's inequality
in time imply
\[
\begin{aligned}
 \norm{e^{-H_a}f}_{L^p}
 &\leq
 \int_{\R}|g(t)|
 \norm{e^{-itH_a}f}_{L^p}\,\dd t\\
 &\leq
 \norm{g}_{L^{r'}}
 \norm{e^{-itH_a}f}_{L_t^rL_x^p}\\
 &\leq
 C\norm{g}_{L^{r'}}\norm{|D|^sf}_{L^2}
 <\infty.
\end{aligned}
\]
We obtain $e^{-H_a}f\in L^p(\R^3)$, contradicting
\eqref{eq:heat-not-Lp}.

When $p=\infty$, the same conclusion follows by duality.
Indeed, for every $\phi\in L^1(\R^3)\cap L^2(\R^3)$, the
representation \eqref{eq:heat-fourier-representation} gives
\[
\begin{aligned}
 |\langle e^{-H_a}f,\phi\rangle|
 &\leq
 \int_{\R}|g(t)|
 |\langle e^{-itH_a}f,\phi\rangle|\,\dd t\\
 &\leq
 \norm{g}_{L^{r'}}
 \norm{e^{-itH_a}f}_{L_t^rL_x^\infty}
 \norm{\phi}_{L^1}\\
 &\leq
 C\norm{g}_{L^{r'}}
 \norm{|D|^sf}_{L^2}\norm{\phi}_{L^1},
\end{aligned}
\]
which implies $e^{-H_a}f\in L^\infty(\R^3)$, contradicting
\eqref{eq:heat-not-Lp}. 
\end{proof}

\noindent{\bf AI disclosure.}
The authors acknowledge the use of AI tools during the exploratory stage
of this project. All mathematical arguments and proofs in the final manuscript were checked and
written by the authors.

\end{document}